\documentclass[12pt]{amsproc}
\usepackage{amsfonts}

\def\M{{\mathcal{M}}}

\def\oM{{{\overline{\mathcal{M}}}}}

\def\CP1{{\mathbb{C}\mathrm{P}^1}}

\def\ll{\langle\langle}

\def\cp{\mathsf{cp\,}}
\def\ml{\mathsf{ml\,}}
\def\nm{\mathsf{nm\,}}
\def\st{\mathsf{st\,}}
\def\ll{\mathsf{ll\,}}

\newtheorem{theorem}{Theorem}[section]
\newtheorem{corollary}[theorem]{Corollary}
\newtheorem{lemma}[theorem]{Lemma}

\numberwithin{equation}{section}

\title{Combinatorics of binomial decompositions of the simplest
Hodge integrals}

\author{S.~V.~Shadrin}

\address{Independent University of Moscow and
Stockholm University}

\email{shadrin@mccme.ru}

\thanks{Partially supported by the grants RFBR 01-01-00660 and
RFBR 02-01-22004.}

\subjclass{14H10}

\date{\today}

\keywords{moduli space of curves, Hodge integrals}

\begin{document}

\begin{abstract}
We reduce the calculation of the simplest Hodge integrals to
some sums over decorated trees. Since Hodge
integrals are already calculated, this gives a proof of a rather
interesting combinatorial theorem and a new representation of
Bernoulli numbers.
\end{abstract}

\maketitle

%\tableofcontents

%%%%%%%%%%%%%%%%%%%%%%%%%%%%%%%%%%%%%%%%%%%%%%%%%%%%%%%%%%%%

%%%%%%%%%%%%%%%%%%%%%%%%%%%%%%%%%%%%%%%%%%%%%%%%%%%%%%%%%%%%

%%%%%%%%%%%%%%%%%%%%%%%%%%%%%%%%%%%%%%%%%%%%%%%%%%%%%%%%%%%%

\section{Introduction}

In this paper we study the simplest Hodge integrals on the
moduli space of curves. Let $\oM_{g,n}$ be the moduli space of
curves of genus $g$ with $n$ marked points. By $\psi_i$ denote
the first Chern class of the line bundle over $\oM_{g,n}$, whose
fiber at a moduli point $(C_g,x_1,\dots,x_n)\in\oM_{g,n}$ is
equal to $T^*_{x_i}C$. The Hodge bundle is the rank $g$ vector
bundle over $\oM_{g,n}$, whose fiber at a moduli point
$(C_g,x_1,\dots,x_n)\in\oM_{g,n}$ is equal to $H^0(C,\omega_C)$.
By $\lambda_1,\dots,\lambda_g$ denote the Chern classes of the
Hodge bundle.

In this paper we consider the integrals
\begin{equation}
\int_{\oM_{g,1}}\psi_1^{3g-2-i}\lambda_i.
\end{equation}
We give an algorithm to calculate these integrals. In fact,
these integrals have already been calculated by Faber and
Pandharipande, see~\cite{fp1}. There is a
remarkable formula for the generating function of Hodge
integrals:
\begin{equation}\label{FPformula}
1+\sum_{g\geq 1}\sum_{i=0}^g t^{2g} k^i
\int_{\oM_{g,1}} \psi_1^{2g-2+i} \lambda_{g-i}
=
\left(
\frac{t/2}{\sin(t/2)}
\right)^{k+1}
.
\end{equation}

So, let us explain our motivations. First, for the integrals
$\int_{\oM_{g,1}} \psi_1^{2g-2} \lambda_{g}$ our computing
algorithm has a rather simple interpretation as a certain sum
over trees. So, since these integrals are already calculated, we
solve an interesting combinatorial problem. In particular, this
gives a representation of Bernoulli numbers as some sums over
trees.

Second, the combinatorial structure of Hodge integrals
still seems to produce some questions, see~\cite{gjv}. Our
approach could clarify something in this direction. Moreover,
the soul of our computing algorithm is a `cut-and-join' type
equation for the Hodge integrals over two-pointed ramification
cycles. Since Hurwitz numbers satisfy the similiar equations,
this could give some new relations between Hurwitz numbers and
Hodge integrals, see~\cite{elsv, gv, llz}.

Third, our computing algorithm is a modification of the
computing algorithm for the simplest Witten intersections,
see~\cite{s}. It works very good to calculate any concrete
intersection number, but it is very hard to prove any general
statement using this approach. In the case considered in this
paper we reduce the computing algorithm to rather friendly
combinatorics. We hope this could help in our approach to the
Witten conjecture.

\subsection*{Organization of the paper}

In Section 2 we give our algorithm for computing the integrals
$\int_{\oM_{g,1}} \psi_1^{2g-2} \lambda_{g}$. In Section 3
explain the combinatorial interpretation of this algorithm
and combinatorial corollaries of this interpretation.
In Section 4 we prove all theorems of Section 2.

In Appendix A we give some calculations checking
independently our algorithm. In Appendix B we describe
separately genus zero case of our combinatorial results.
In Appendix C we generalize our algorithm to compute all
Hodge integrals $\int_{\oM_{g,1}} \psi_1^{3g-2-i}
\lambda_{i}$.

%%%%%%%%%%%%%%%%%%%%%%%%%%%%%%%%%%%%%%%%%%%%%%%%%%%%%%%%%%%%

%%%%%%%%%%%%%%%%%%%%%%%%%%%%%%%%%%%%%%%%%%%%%%%%%%%%%%%%%%%%

%%%%%%%%%%%%%%%%%%%%%%%%%%%%%%%%%%%%%%%%%%%%%%%%%%%%%%%%%%%%

\section{Calculation of Hodge integrals}\label{algorithm}

In this section we explain an algorithm to calculate the Hodge
integrals $\int_{\oM_{g,1}}\psi_1^{2g-2}\lambda_g$. The proofs
of all theorems of this section are given in Section 4.

\subsection{Some notations}\label{sn}
Consider the moduli space of curves $\oM_{g,n+1}$. Let us fix
some positive integers $a_1,\dots,a_n$. By $V^\circ_g
(a_1,\dots,a_n)$ denote the subvariety of the open moduli space
$\M_{g,n+1}$ consisting of curves $(C,x_1,\dots,x_{n+1})$ such
that
%\begin{equation}
$-(\sum_{i=1}^n a_i) x_1 + a_1 x_2 + \dots + a_n x_{n+1}$
%\end{equation}
is a divisor of a meromorphic function.

Let $V_g (a_1,\dots,a_n)$ be the closure of $V^\circ_g
(a_1,\dots,a_n)$ in $\oM_{g,n}$. By $W_g(\prod_{i=1}^n
\eta_{a_i})$ denote the integral of $\psi_1^{g+n-2}\lambda_g$
over the subspace $V_g (a_1,\dots,a_n)$:

\begin{equation}
W_g(\prod_{i=1}^n \eta_{a_i})\colon =
\int_{V_g (a_1,\dots,a_n)}\psi_1^{g+n-2}\lambda_g.
\end{equation}

\subsection{Binomial decomposition}

\begin{theorem}\label{binomial}
For arbitrary positive integers
$a_1,\dots,a_n$ we have
\begin{align}
(-1)^g g! \int_{\oM_{g,1}}\psi_1^{2g-2}\lambda_g = &
\binom{g}{0}W_g(\prod_{i=1}^n \eta_{a_i})
- \binom{g}{1}
W_g(\eta_1\prod_{i=1}^n \eta_{a_i})
\\
& + \dots
+ (-1)^g \binom{g}{g}
W_g(\eta_1^g\prod_{i=1}^n \eta_{a_i}).\notag
\end{align}
\end{theorem}

\subsection{Recursion relation}
There is a recursion relation for the
the numbers $W_g(\prod_{i=1}^n \eta_{a_i})$:

\begin{theorem}\label{recursion}
If $g+n-2>0$, then
\begin{align}\label{rec}
(\sum_{i=1}^n a_i)(2g+n-1) W_g(\prod_{i=1}^n \eta_{a_i})=
&
\sum_{k<l} (a_k+a_l) W_g(\eta_{a_k+a_l}\prod_{i\not= k,l}
 \eta_{a_i}) \\
& +
\sum_{k=1}^n \frac{a_k^3-a_k}{12}
W_{g-1}(\prod_{i=1}^n \eta_{a_i}). \notag
\end{align}
\end{theorem}

\subsection{Initial values}

Consider the case $g+n-2=0$. This means that either $g=1$,
$n=1$, or $g=0$, $n=2$. In these cases we have the following:

\begin{theorem}\label{initial}
\begin{equation}
W_1(\eta_{a_1})=\frac{a_1^2-1}{24}; \qquad
W_0(\eta_{a_1}\eta_{a_2})=1.
\end{equation}
\end{theorem}

\subsection{Calculation of Hodge integrals}

Now we have an algorithm to calculate the integrals
$\int_{\oM_{g,1}}\psi_1^{2g-2}\lambda_g$. First, we
express such integral via the numbers $W_g(\prod_{i=1}^n
\eta_{a_i})$ (Theorem~\ref{binomial}). Then we step by step
simplify the numbers $W_g(\prod_{i=1}^n
\eta_{a_i})$ using Theorem~\ref{recursion} until we obtain
numbers calculated in Theorem~\ref{initial}.

Let us apply this algorithm to calculate
$\int_{\oM_{1,1}}\lambda_1$.
The first step looks, for example, like follows:
$$
-\int_{\oM_{1,1}}\lambda_1=W_1(\eta_2)-W_1(\eta_1\eta_2).
$$
Then
$$
3\cdot 3\cdot W_1(\eta_1 \eta_2)= 3 \cdot W_1(\eta_3) +
 \frac{8-2}{12}
W_0(\eta_1 \eta_2).
$$
Since
$$
W_0(\eta_1 \eta_2)=1, \quad W_1(\eta_3)=\frac{8}{24}, \quad
\mathrm{and} \quad W_1(\eta_2)=\frac{3}{24},
$$
we have
$$
-\int_{\oM_{1,1}}\lambda_1=
\frac{3}{24}-\frac{1}{3}\cdot\frac{8}{24}-\frac{1}{9}\cdot
\frac{6}{12} =-\frac{1}{24}.
$$

The same number is give by Equation~\ref{FPformula}.

\section{The number $W_g(\eta_1^n)$ as a sum over trees}

The algorithm of calculation of numbers $W_g(\eta_1^n)$
implies a representation of these numbers as a sum over
some trees. In this section we describe the trees we need and
prove a theorem expressing $W_g(\eta_1^n)$ via some numbers
calculated by trees.

\subsection{Decorated trees}
Let us fix $g\geq 0$ and $n\geq 1$. We describe here trees
corresponding to the number $W_g(\eta_1^n)$.

\subsubsection{Construction}
We consider rooted trees. Each vertex has $\leq 2$ sons.
\footnote{If we have two connected vetices in a rooted tree,
then one of these vertices is the son of another one.
By son we call the vertex, which is farther from the root.}
We require that there are exactly $n$ vertices with no
sons and exactly $g$ vertices with one son. This follows that
there are exactly $n-1$ vertices with two sons.

By leaves we call vertices with no sons. By $V_0$ denote the set
of leaves. Let us assign to each leaf a personal number $\nm$
from the set $\{1,\dots,n\}$. In other words, we consider a
one-to-one correspondence $\nm\colon V_0 \to \{1,\dots,n\}$.

One more condition for the structure of trees is that the father
of each leaf has two sons (In fact, we use this condition just
to decrease the number of examples. The contribution of such
trees to our formulas will be zero).

By $V_1$ denote the set of vertices with one son.
By $V_2$ denote the set of vertices with two sons.
Consider a map $\cp\colon V_1\cup V_2\to\{1,\dots,2g+n-1\}$
satifying the following three properties:

\begin{enumerate}
\item The map $\cp$ takes different verices to different
numbers.

\item If $a\in \cp(V_1)$, then $a>1$ and $a-1\not\in \cp(V_1\cup
V_2)$.

\item If vertex $v$ is a `descendant' of vertex $v'$, then
$\cp(v)>\cp(v')$.
\end{enumerate}

In particular, if the root vertex has one (resp., two) sons,
then its image under the mapping $\cp$ is equal to $2$ (resp.,
$1$).

By a decorated tree (or $(n,g)$-decorated tree) we call all this
data, i.~e. a rooted tree with mappings $\nm$ and $\cp$.

\subsubsection{Examples}

Consider $n=3$, $g=2$. Up to isomorphism, there are $3$
possible rooted trees with the mapping $\cp$.

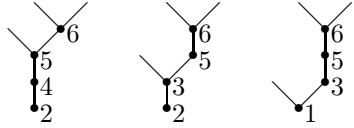
\begin{picture}(150,60)
\put(20,10){\line(0,1){20}}
\put(20,30){\line(-1,1){10}}
\put(20,30){\line(1,1){20}}
\put(30,40){\line(-1,1){10}}
\put(20,10){\circle*{3}}
\put(22,5){\footnotesize 2}
\put(20,20){\circle*{3}}
\put(22,15){\footnotesize 4}
\put(20,30){\circle*{3}}
\put(22,25){\footnotesize 5}
\put(30,40){\circle*{3}}
\put(32,35){\footnotesize 6}
\put(70,10){\line(0,1){10}}
\put(70,20){\line(-1,1){10}}
\put(70,20){\line(1,1){10}}
\put(80,30){\line(0,1){10}}
\put(80,40){\line(-1,1){10}}
\put(80,40){\line(1,1){10}}
\put(70,10){\circle*{3}}
\put(72,5){\footnotesize 2}
\put(70,20){\circle*{3}}
\put(72,15){\footnotesize 3}
\put(80,30){\circle*{3}}
\put(82,25){\footnotesize 5}
\put(80,40){\circle*{3}}
\put(82,35){\footnotesize 6}
\put(120,10){\line(-1,1){10}}
\put(120,10){\line(1,1){10}}
\put(130,20){\line(0,1){20}}
\put(130,40){\line(-1,1){10}}
\put(130,40){\line(1,1){10}}
\put(120,10){\circle*{3}}
\put(122,5){\footnotesize 1}
\put(130,20){\circle*{3}}
\put(132,15){\footnotesize 3}
\put(130,30){\circle*{3}}
\put(132,25){\footnotesize 5}
\put(130,40){\circle*{3}}
\put(132,35){\footnotesize 6}
\end{picture}

Note that in this particular case the type of rooted graph
uniquely determines the mapping $\cp$. For each of these graphs
there are three possible mappings $\nm$.

Consider $n=2$, $g=2$. There is a unique possible rooted tree
with the mapping $\cp$. For this tree there is a unique
possible mapping $\nm$.

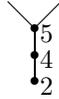
\begin{picture}(150,50)
\put(20,10){\line(0,1){20}}
\put(20,30){\line(-1,1){10}}
\put(20,30){\line(1,1){10}}
\put(20,10){\circle*{3}}
\put(22,5){\footnotesize 2}
\put(20,20){\circle*{3}}
\put(22,15){\footnotesize 4}
\put(20,30){\circle*{3}}
\put(22,25){\footnotesize 5}
\end{picture}

Consider $n=1$, $g=2$. There exists no such decorated trees.
Moreover, if $n=1$, then a decorated tree exists if and only if
$g=0$. In this case this tree consists just of one vertex.

\subsection{Taking a number from a decorated tree}

\subsubsection{Function $ml$.}
Consider a decorated tree. Let us define one more function
$\ml\colon V_1\cup V_2\to\{1,\dots,n\}$. The mapping $\ml$ takes
a vertex to the number of its descendants in $V_0$.

In the examples studied above the mapping $\ml$ looks like
follows:

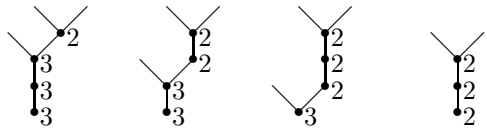
\begin{picture}(200,60)
\put(20,10){\line(0,1){20}}
\put(20,30){\line(-1,1){10}}
\put(20,30){\line(1,1){20}}
\put(30,40){\line(-1,1){10}}
\put(20,10){\circle*{3}}
\put(22,5){\footnotesize 3}
\put(20,20){\circle*{3}}
\put(22,15){\footnotesize 3}
\put(20,30){\circle*{3}}
\put(22,25){\footnotesize 3}
\put(30,40){\circle*{3}}
\put(32,35){\footnotesize 2}
\put(70,10){\line(0,1){10}}
\put(70,20){\line(-1,1){10}}
\put(70,20){\line(1,1){10}}
\put(80,30){\line(0,1){10}}
\put(80,40){\line(-1,1){10}}
\put(80,40){\line(1,1){10}}
\put(70,10){\circle*{3}}
\put(72,5){\footnotesize 3}
\put(70,20){\circle*{3}}
\put(72,15){\footnotesize 3}
\put(80,30){\circle*{3}}
\put(82,25){\footnotesize 2}
\put(80,40){\circle*{3}}
\put(82,35){\footnotesize 2}
\put(120,10){\line(-1,1){10}}
\put(120,10){\line(1,1){10}}
\put(130,20){\line(0,1){20}}
\put(130,40){\line(-1,1){10}}
\put(130,40){\line(1,1){10}}
\put(120,10){\circle*{3}}
\put(122,5){\footnotesize 3}
\put(130,20){\circle*{3}}
\put(132,15){\footnotesize 2}
\put(130,30){\circle*{3}}
\put(132,25){\footnotesize 2}
\put(130,40){\circle*{3}}
\put(132,35){\footnotesize 2}
\put(180,10){\line(0,1){20}}
\put(180,30){\line(-1,1){10}}
\put(180,30){\line(1,1){10}}
\put(180,10){\circle*{3}}
\put(182,5){\footnotesize 2}
\put(180,20){\circle*{3}}
\put(182,15){\footnotesize 2}
\put(180,30){\circle*{3}}
\put(182,25){\footnotesize 2}
\end{picture}

\subsubsection{The numbers $N(\Gamma)$ and $S_{g,n}$.}

To any $(n,g)$-decorated tree $\Gamma$ we assign the number
$N(\Gamma)$:

\begin{equation}
N(\Gamma)=\frac{1}{n^{(n+g-1)}}\cdot \prod_{v\in V_2}
\frac{\ml(v)}{\cp(v)} \cdot \prod_{v\in V_1}
\frac{\ml(v)^3-\ml(v)}{12\cp(v)}
\end{equation}

By $S_{g,n}$ denote the sum
\begin{equation}
S_{g,n}=\sum_{\Gamma} N(\Gamma),
\end{equation}
where the sum is taken over all $(n,g)$-decorated trees
$\Gamma$.

\subsection{Theorems}

\begin{theorem}\label{WeqS}
For any $n\geq 1$, $g\geq 0$, we have
\begin{equation}
W_g(\eta_1^n)=S_{g,n}.
\end{equation}
\end{theorem}

This theorem implies purely combinatorial identities.

\begin{corollary} \label{arbgenus}
For any $g\geq 1$, $n\geq 1$, we have
\begin{multline}
\binom{g}{0}S_{g,n+g}
- \binom{g}{1}
S_{g,n+g-1}
+
\dots
+ (-1)^g \binom{g}{g}
S_{g,n}=\\
\frac{(2^{2g-1}-1) g!}
{2^{2g-1} (2g)!}
|B_{2g}|.
\end{multline}
\end{corollary}

\begin{corollary} \label{genuszero}
For any $n\geq 1$, we have
\begin{equation}
S_{0,n}=1.
\end{equation}
\end{corollary}

Note that the trees determining the numbers $S_{0,n}$ have much
more simple description then in the general case. So
Corollary~\ref{genuszero} has more simple and more beautiful
reformulation.

\subsection{Proof of Theorem~\ref{WeqS}}

In fact, this Theorem is more or less obvious.
Let us consider the calculation of the number $W_g{\eta_1^n}$ step by step.
Simultaniously we try to build the corresponding graph.

At the beginning we have only $n$ vertices $v_1^0,\dots,v_n^0$.
The first step given by Equation~\ref{recursion}:
\begin{equation}
W_g{\eta_1^n}=\frac{1}{n(2g+n-1)}\sum_{i<j}2W_g(\eta_2\eta_1^{n-2}).
\end{equation}
Consider the summand in the right hand side of
this equation corresponding to the pair $(i,j)$, $i<j$.
According to this summand we add to our graph
the vertex $v_{2g+n-1}^2$ with two sons, $v_i^0$ and $v_j^0$.
We put $\cp(v_{2g+n-1}^2)=2g+n-1$ and we have $\ml(v_{2g+n-1}^2)=2$.
Therefore, the factor $2/(2g+n-1)$ is equal to
$\ml(v_{2g+n-1}^2)/\cp(v_{2g+n-1}^2)$
(Let us always skip the factor $1/n$ throughout our calculations).

Let us make the second step for this summand. We have:
\begin{align}
W_g(\eta_2\eta_1^{n-2})=&\frac{1}{n(2g+n-2)}\sum_{i'<j'}2W_g(\eta_2^2\eta_1^{n-4})\\
&+\frac{1}{n(2g+n-2)}\sum_{i'}3W_g(\eta_3\eta_1^{n-3})\notag\\
&+\frac{1}{n(2g+n-2)}\cdot
\frac{2^3-2}{12}\cdot
W_{g-1}(\eta_2\eta_1^{n-2}).\notag
\end{align}
Here we take all sums over $i',j'\in\{1,\dots,n\}\setminus\{i,j\}$.

Consider the summand of the first sum in the right hand side of
this equation corresponding to the pair $(i',j')$, $i'<j'$.
According to this summand we add to our graph
the vertex $v_{2g+n-2}^2$ with two sons, $v_{i'}^0$ and $v_{j'}^0$.
We put $\cp(v_{2g+n-2}^2)=2g+n-2$ and we have $\ml(v_{2g+n-2}^2)=2$.
Therefore, the factor $2/(2g+n-2)$ is equal to
$\ml(v_{2g+n-2}^2)/\cp(v_{2g+n-2}^2)$.

Consider the summand of the second sum in the right hand side of
this equation corresponding to $i'$.
According to this summand we add to our graph
the vertex $v_{2g+n-2}^2$ with two sons, $v_{i'}^0$ and $v_{2g+n-1}^2$.
We put $\cp(v_{2g+n-2}^2)=2g+n-2$ and we have $\ml(v_{2g+n-2}^2)=3$.
Therefore, the factor $3/(2g+n-2)$ is equal to
$\ml(v_{2g+n-2}^2)/\cp(v_{2g+n-2}^2)$.

Consider the summand of the second sum in the right hand side of
this equation corresponding to $W_{g-1}(\eta_2\eta_1^{n-2})$.
According to this summand we add to our graph
the vertex $v_{2g+n-2}^1$ with one sons, $v_{2g+n-1}^2$.
We put $\cp(v_{2g+n-2}^1)=2g+n-2$ and we have $\ml(v_{2g+n-2}^1)=2$.
Therefore, the factor $(2^3-2)/12(2g+n-2)$ is equal to
$\ml(v_{2g+n-2}^1)/\cp(v_{2g+n-2}^1)$.

Continue with this procedure one obtain the same representation
of $W_g(\eta_1^n)$ as we have described. The missed factor $1/n$ gives
the contribution $1/n^{g+n-1}$; the same factor we have in our definition of
$N(\Gamma)$.

\section{Proofs of Theorems~\ref{binomial}-\ref{initial}}

\subsection{Proof of Theorem~\ref{initial}}

This Theorem consists of two formulas:
\begin{equation}
W_1(\eta_{a_1})=\frac{a_1^2-1}{24}; \qquad
W_0(\eta_{a_1}\eta_{a_2})=1.
\end{equation}

Let us start with the second one. We have
$W_0(\eta_{a_1}\eta_{a_2})=\int_{V_0 (a_1,a_2)}1$.
Since $V_0(a_1,a_2)=\oM_{0,3}$, it follows the required
formula.

Consider now the first formula. It could be rewritten as
\begin{equation}
\int_{V_1(a_1)}\lambda_1=\frac{a_1^2-1}{24}.
\end{equation}
Note that $\lambda_1|_{V_1(a_1)}=\psi_1|_{V_1(a_1)}$.
The integral $\int_{V_1(a_1)}\psi_1$ equals
\begin{equation}
\frac{1}{4a_1}\sum_{i=1}^{a_1-1}i(a_1-i)\int_{V_0(i,a_1-i)}1
\end{equation}
(this follows from~\cite{s}, Theorem 12.2).
Since $W_0(i,a_1-i)=\int_{V_0(i,a_1-i)}1=1$ and
$\sum_{i=1}^{a_1-1}i(a_1-i)=(a_1^3-a_1)/6$, it follows that
\begin{equation}
\int_{V_1(a_1)}\lambda_1=\int_{V_1(a_1)}\psi_1=
\frac{a_1^2-1}{24}.
\end{equation}

\subsection{Ionel Lemma}

In this subsection we recall the Ionel Lemma, which is the main
tool in the foregoing argument. For a more detailed explanation
of this technique, see~\cite{i,s}.

Consider the space $\widehat H$ of admissible
coverings of genus $g$, of degree $n$, with $m$ critical values,
and with partitions $A_1,\dots,A_m$,
$A_i=(a^i_1,\dots,a^i_{l_i})$, over these critical values. Of
course, for any $i$, $\sum_{j=1}^{l_i}a^i_j=n$, and
$\sum_{i=1}^m\sum_{j=1}^{l_i}(a^i_j-1)=
nm-\sum_{i=1}^ml_i=2g+2n-2$. Moreover, we consider admissible
coverings with all marked preimages of all critical values.

There are two mapping of the space $\widehat H$. The first one,
$\st\colon \widehat H\to \oM_{g,\sum_{i=1}^m l_i}$, takes an
admissible covering to its source curve with marked preimages of
critical values. The second mapping, $\ll\colon \widehat H\to
\oM_{0,m}$, takes an admissible covering to its target curve
with marked critical values.

Consider an admissible covering $f\in\widehat H$. Let
$$
\st(f)=(C_g,x^1_1,\dots,x^1_{l_1},\dots,
x^m_1,\dots,x^m_{l_m}), \quad \ll(f)=(C_0,z_1,\dots,z_m).
$$
We choose our notations to make $x^i_j$ be a preimage of
$x_i$ of multiplicity $a^i_j$ w.~r.~t. the covering $f$.

Let $\psi(x^i_j)$ be the first Chern class of the line bundle
over $\oM_{g,\sum_{i=1}^m l_i}$, whose fiber at a moduli point
$(C_g,x^1_1,\dots,x^m_{l_m})$ is equal to $T^*_{x^i_j}C_g$.
Let $\psi(z_i)$ be the first Chern class of the line bundle
over $\oM_{o,m}$, whose fiber at a moduli point
$(C_0,z_1,\dots,z_m)$ is equal to $T^*_{z_i}C_g$.

\begin{lemma} {\rm (Ionel Lemma, \cite{i})}
In cohomology ring of $\widehat H$ we have:
\begin{equation}
a^i_j\st^*\psi(x^i_j)=\ll^*\psi(z_i).
\end{equation}
\end{lemma}

\subsection{Proof of Theorem~\ref{recursion}}

\subsubsection{}
Consider the space $V_g(a_1,\dots,a_n)$. Let $\widehat H$ be the
space of admissible coverings of genus $g$, of degree
$N=a_1+\dots+a_n$, with $m=2g+n+1$ critical values, and with
partitions $A_1=(N)$, $A_2=(a_1,\dots,a_n)$,
$A_3=\dots=A_m=(2,1,\dots,1)$ over these critical values.

In the foregoing we use the notations from the previous
subsection.

Consider the projection $\pi\colon \oM_{g, \sum_{i=1}^m l_i}\to
\oM_{g, 1+n}$, forgetting all marked points except for $x^1_1,
x^2_1,\dots,x^2_n$. Note that $\pi\circ\st(\widehat
H)=V_g(a_1,\dots,a_n)$. Moreover, $\pi_*\st_*[\widehat H]=
(m-2)![V_g(a_1,\dots,a_n)]$.

We know that $\psi(x^1_1)|_{\st(\widehat H)}=\pi^*\psi(x^1_1)$
(this is proved in~\cite{s}). Therefore,
\begin{equation}
(m-2)!\psi(x_1^1)|_{V_g(a_1,\dots,a_n)}=
\pi_*\st_*\st^*\psi(x^1_1).
\end{equation}
Since
$N\st^*\pi^*\psi(x^1_1)=\ll^*\psi(z_1)$ (this is Ionel Lemma),
it follows that
\begin{equation}
(m-2)!\psi(x_1^1)|_{V_g(a_1,\dots,a_n)}=
\frac{1}{N}\pi_*\st_*\ll^*\psi(z_1).
\end{equation}
So, we are to take a divisor dual to $\psi(z_1)$ on $\oM_{0,m}$,
then we get its preimage in $\widehat H$, and the mapping
$\pi\circ\st$ takes this preimage to the divisor in
$V_g(a_1,\dots,a_n)$ dual to
$\psi(x_1^1)|_{V_g(a_1,\dots,a_n)}$.

Note that $\psi(z_1)$ on $\oM_{0,m}$ is dual to the divisor,
whose generic point is represented by a two-component curve
such that $z_1$ lie on one component and $z_2$ and $z_3$ lie on
the other component.

\subsubsection{}
Consider $\psi(x^1_1)^{g+k-2}\lambda_g\cdot
[V_g(a_1,\dots,a_k)]$. If $g+k-2=0$, then we have one of the
cases considered in Theorem~\ref{initial}. Suppose that
$g+k-2>0$. Then we have
\begin{multline}
\psi(x^1_1)^{g+k-2}\lambda_g\cdot
[V_g(a_1,\dots,a_k)]= \\
\frac{1}{N (m-2)!}
\psi(x^1_1)^{g+k-3}\lambda_g\cdot
\pi_*\st_*\ll^*(\psi(z_1)\cdot [\oM_{0,m}])
\end{multline}

From dimensional conditions and since $\lambda_g$ restricted to
the divisor of irreducible self-intersecting curves equals zero,
it follows that only two types of divisors contribute to
$\psi(x^1_1)^{g+k-3}\lambda_g\cdot
\pi_*\st_*\ll^*(\psi(z_1)\cdot \oM_{0,m})$.
(We skip here some argument on the boundary behaviour of the
space $\widehat H$. We refer to~\cite{s} for details.)

A divisor $D_{i,j}$ of the first type consists of two-compenent
curves $(C,x^1_1,x^2_1,\dots,x^2_n)\in V_g(a_1,\dots,a_n)$ such
that one component of $C$ has genus zero and contains
only $x^2_i$ and $x^2_j$. Here pairs $(i,j)$ enumerate
such divisors. Obviously, the restriction of
$\psi(x^1_1)^{g+k-3}\lambda_g$ to the divisor $D_{i,j}$ is equal
to the integral of this class over
$$V_g(a_1,\dots,\hat a_i,\dots,\hat a_j,\dots, a_n, a_i+a_j).$$

The multiplicity of the mapping $\pi\circ\st$ over $D_{i,j}$
equals $(m-3)!$. The multiplicity of the mapping $\ll$ at
$(\pi\circ\st)^{-1}(D_{i,j})$ equals $a_i+a_j$. Thus we obtain
that the coefficient of the correspoding summand in
Equation~\ref{rec} equals
\begin{equation}
\frac{(a_i+a_j)(m-3)!}{N(m-2)!}
=
\frac{(a_i+a_j)}{(2g+n-1)(\sum_{i=1}^n a_i)}.
\end{equation}

A divisor $D_i$ of the second type consists of two-component
curves $(C,x^1_1,x^2_1,\dots,x^2_n)\in V_g(a_1,\dots,a_n)$
such that one component of $C$ has genus one and contains
only $x_i$. Obviously, the restriction of
$\psi(x^1_1)^{g+k-3}\lambda_g$ to the divisor $D_{i,j}$ is equal
to the integral of $\psi(x^1_1)^{g+k-3}\lambda_{g-1}$
over $V_{g-1}(a_1,\dots,a_n)$ multiplied by the integral
of $\lambda_1$ over $V_1(a_i)$. It follows from
Theorem~\ref{initial} that last factor equals $(a_i^2-1)/24$.

The multiplicity of the mapping $\pi\circ\st$ over $D_{i}$
equals $2\cdot (m-3)!$. The multiplicity of the mapping $\ll$ at
$(\pi\circ\st)^{-1}(D_{i,j})$ equals $a_i$. Thus we
obtain that the coefficient of the correspoding summand in
Equation~\ref{rec} equals
\begin{equation}
a_i\cdot \frac{(a_i^2-1)}{24}\cdot \frac{2(m-3)!}{N(m-2)!}
=
\frac{a_i^3-a_i}{12(2g+n-1)
(\sum_{i=1}^n a_i)}.
\end{equation}

This concludes the proof.

\subsection{Proof of Theorem~\ref{binomial}}

\subsubsection{}
Let us fix $g$ and $a_1,\dots,a_n$. By $V_g^i$ denote
\begin{equation}
V_g^i :=
V_g(a_1,\dots,a_n\underbrace{1,\dots,1}_i).
\end{equation}
Recall that $V_g^i$ is the subspace of
$\oM_{g,n+i+1}$ consisting of curves
$(C,x_1,\dots,x_{n+i+1})$ such that
$$
-(N+i)x_1+a_1x_2+\dots+a_nx_{n+1}+x_{n+2}+\dots+x_{n+i+1}
$$
is the divisor of a meromorphic function.

Consider the mappings $\pi_{j,i}\colon V_g^j\to\oM_{g,n+j+1-i}$
forgetting the marked points $x_{n+j+2-i},\dots,x_{n+j+1}$.
Note that $\pi_{j,i,*}[V^j_g]=i![\pi_{i,*}(V^j_g)]$,
$\pi_{j,i}(V^j_g)=\pi_{j+1,i+1}(V^{j+1}_g)$,
and $\pi_{g,g}(V_g^g)=\oM_{g,n+1}$.

Recall that our goal is to get an expression of
$\int_{\oM_{g,1}}\psi_1^{2g-2}\lambda_g=
\int_{\oM_{g,n+1}}\psi_1^{2g+n-2}\lambda_g$.
We shall do this in $g$ steps.

\subsubsection{}
Consider the projection $\sigma\oM_{g,n+2}\to\oM_{g,n+1}$
forgetting the last marked point. We have
\begin{equation}
\sigma^*(\psi_1)^{2g+n-2}=\psi_1^{2g+n-2}-
\sigma^*(\psi_1)^{2g+n-3}\cdot D,
\end{equation}
where $D$ is the class of the divisor, whose generic point
is represented by a two-component curve such that one component
has genus zero and contains $x_1$ and $x_{n+2}$ and the other
component has genus $g$ and contains all other points.

Since $\sigma^*\lambda=\lambda$, we have:
\begin{multline}
I:=\int_{\oM_{g,n+1}}\psi_1^{2g+n-2}\lambda_g=\\
\int_{\pi_{g,g}(V_g^g)}\psi_1^{2g+n-2}\lambda_g=
\frac{1}{g}\int_{\pi_{g,g-1}(V_g^g)}\sigma^*(\psi_1)^{2g+n-2}
\lambda_g
\end{multline}

Note that
\begin{equation}
\sigma(D\cap\pi_{g,g-1}(V_g^g))
=\pi_{g-1,g-1}(V_g^{g-1}).
\end{equation}
Therefore,
\begin{equation}\label{expI}
I=\frac{1}{g}\int_{\pi_{g,g-1}(V_g^g)}\psi_1^{2g+n-2}
\lambda_g-
\frac{1}{g}\int_{\pi_{g-1,g-1}(V_g^{g-1})}\psi_1^{2g+n-3}
\lambda_g.
\end{equation}

\subsubsection{}
Applying the same argument to the right hand
side of Equation~\ref{expI} we get:
\begin{align}
\int_{\pi_{g,g-1}(V_g^g)}\psi_1^{2g+n-2}
\lambda_g= &
\frac{1}{g-1}\int_{\pi_{g,g-2}(V_g^g)}\psi_1^{2g+n-2}
\lambda_g\\
&-\frac{1}{g-1}
\int_{\pi_{g-1,g-2}(V_g^{g-1})}\psi_1^{2g+n-3}
\lambda_g\notag
\end{align}
and
\begin{align}
\int_{\pi_{g-1,g-1}(V_g^{g-1})}\psi_1^{2g+n-3}
\lambda_g=&
\frac{1}{g-1}\int_{\pi_{g-1,g-2}(V_g^{g-1})}\psi_1^{2g+n-3}
\lambda_g\\
&-\frac{1}{g-1}
\int_{\pi_{g-2,g-2}(V_g^{g-2})}\psi_1^{2g+n-4}
\lambda_g\notag
\end{align}

Thus we have:
\begin{align}
g(g-1)I=&\int_{\pi_{g,g-2}(V_g^g)}
\psi_1^{2g+n-2}\lambda_g
-2\int_{\pi_{g-1,g-2}(V_g^{g-1})}\psi_1^{2g+n-3}
\lambda_g\\
&+\int_{\pi_{g-2,g-2}(V_g^{g-2})}\psi_1^{2g+n-4}
\lambda_g.\notag
\end{align}

Continue with this procedure we obtain:
\begin{align}
g!I=&\binom{g}{0}\int_{\pi_{g,0}(V_g^g)}
\psi_1^{2g+n-2}\lambda_g
-\binom{g}{1}
\int_{\pi_{g-1,0}(V_g^{g-1})}\psi_1^{2g+n-3}
\lambda_g \\
&+\dots
+(-1)^g\binom{g}{g}
\int_{\pi_{0,0}(V_g^{0})}\psi_1^{g+n-2}
\lambda_g.\notag
\end{align}
Since $\pi_{i,0}$ are identical mappings, this equation is
exactly the statement of Theorem~\ref{binomial}.

\appendix

\section{Some calculations in genus 2}

Here we give some calculations checking independently in a
special case Corollary~\ref{arbgenus} and therefore the whole
algorithm for computing Hodge integrals.

Let us put $g=2$, $n=1$. Recall that $B_{4}=-1/30$. So,
we would like to check that
\begin{equation}
S_{2,3}-2 S_{2,2}=
\frac{7\cdot 2}
{8\cdot 24\cdot 30}
=
\frac{7}{2^6\cdot 3^2\cdot 5}
\end{equation}

We have already studied in examples all types of graphs
contributing to $S_{2,3}$ and $S_{2,2}$:

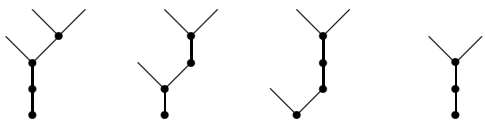
\begin{picture}(200,60)
\put(20,10){\line(0,1){20}}
\put(20,30){\line(-1,1){10}}
\put(20,30){\line(1,1){20}}
\put(30,40){\line(-1,1){10}}
\put(20,10){\circle*{3}}
\put(20,20){\circle*{3}}
\put(20,30){\circle*{3}}
\put(30,40){\circle*{3}}
\put(70,10){\line(0,1){10}}
\put(70,20){\line(-1,1){10}}
\put(70,20){\line(1,1){10}}
\put(80,30){\line(0,1){10}}
\put(80,40){\line(-1,1){10}}
\put(80,40){\line(1,1){10}}
\put(70,10){\circle*{3}}
\put(70,20){\circle*{3}}
\put(80,30){\circle*{3}}
\put(80,40){\circle*{3}}
\put(120,10){\line(-1,1){10}}
\put(120,10){\line(1,1){10}}
\put(130,20){\line(0,1){20}}
\put(130,40){\line(-1,1){10}}
\put(130,40){\line(1,1){10}}
\put(120,10){\circle*{3}}
\put(130,20){\circle*{3}}
\put(130,30){\circle*{3}}
\put(130,40){\circle*{3}}
\put(180,10){\line(0,1){20}}
\put(180,30){\line(-1,1){10}}
\put(180,30){\line(1,1){10}}
\put(180,10){\circle*{3}}
\put(180,20){\circle*{3}}
\put(180,30){\circle*{3}}
\end{picture}

The number $S_{2,3}$ is determined by the first three types of
graphs. These graphs have 3 possible mappings
$\nm$. The mappings $\cp$ and $\ml$ for these graphs are defined
above. Thus we have that the first type of graphs contributes
\begin{equation}
3\cdot \frac{1}{3^4} \frac{2\cdot 3}{2\cdot 4\cdot 5\cdot 6}
\left(\frac{24}{12}\right)^2
=\frac{1}{2\cdot 3^3\cdot 5},
\end{equation}
the second type of graphs contributes
\begin{equation}
3\cdot \frac{1}{3^4} \frac{2\cdot 3}{2\cdot 3\cdot 5\cdot 6}
\frac{6\cdot 24}{12^2}
=\frac{1}{2\cdot 3^4\cdot 5},
\end{equation}
and the third type of graphs contributes
\begin{equation}
3\cdot \frac{1}{3^4} \frac{2\cdot 3}{1\cdot 3\cdot 5\cdot 6}
\left(\frac{6}{12}\right)^2
=\frac{1}{2^2\cdot 3^4\cdot 5}.
\end{equation}

Therefore,
\begin{equation}
S_{2,3}=\frac{1}{2^2\cdot 3^2\cdot 5}.
\end{equation}

The number $S_{2,2}$ is determined by the fourth type of
graph. This graph has the unique possible mapping
$\nm$. The mappings $\cp$ and $\ml$ for this graph are
also defined above. Thus we have
\begin{equation}
S_{2,2}=\frac{1}{2^3}\frac{2}{2\cdot 4\cdot 5}
\left(\frac{6}{12}\right)^2=\frac{1}{2^7\cdot 5}.
\end{equation}

Therefore,
\begin{equation}
S_{2,3}-2S_{2,2}=\frac{1}{2^2\cdot 3^2\cdot 5}-
\frac{1}{2^6\cdot 5}=\frac{7}{2^6\cdot 3^2\cdot 5}.
\end{equation}
The same equality is given by Corollary~\ref{arbgenus}.

\section{Genuz $0$ case}

In genus zero we have the most simple and beautiful
combinatorics. Moreover, in this case our argument with
intersections on the moduli spaces of curves could be
considered as a purely combinatorial argument.

\subsection{Genus zero combinatorial statement}

Let us fix $g\geq 0$ and $n\geq 1$. We describe here trees
corresponding to the number $W_g(\eta_1^n)$.

We consider rooted trees. Each vertex has either $2$ sons or
no sons. We have exactly $n$ vertices with
no sons (leaves) an $n-1$ vertices with two sons.

By $V_0$ denote the set of leaves. An item of the
construction is a one-to-one correspondence $\nm\colon
V_0 \to \{1,\dots,n\}$.

By $V_2$ denote the set of vertices with two sons.
Consider a one-to-one map $\cp\colon V_2\to\{1,\dots,n-1\}$
satifying the following property:
If vertex $v$ is a `descendant' of vertex $v'$, then
$\cp(v)>\cp(v')$.

By a decorated $n$-tree we call all this data, i.~e. a rooted
tree with mappings $\nm$ and $\cp$.

Consider a decorated $n$-tree. Let us define one more function
$\ml\colon V_2\to\{1,\dots,n\}$. The mapping $\ml$ takes
a vertex to the number of its descendants in $V_0$.

Let $\Gamma$ be a decorated tree. By $N(\Gamma)$ denote the
number
\begin{equation}
N(\Gamma):=\frac{1}{n^{(n-1)}}\prod_{v\in V_2}
\frac{\ml(v)}{\cp(v)}
\end{equation}

The special case of Theorem~\ref{WeqS} looks like follows:

\begin{theorem}\label{gZst}
For any $n$ we have
\begin{equation}
\sum_{\Gamma}N(\Gamma)=1,
\end{equation}
where the sum is taken over all decorated $n$-trees $\Gamma$.
\end{theorem}

\subsection{Proof of Theorem~\ref{gZst}}

As we have already shown in Section 3,
$\sum_{\Gamma}N(\Gamma)$ equals the number $W_0(\eta_1^n)$,
where the numbers $W_0(\prod_{i=1}^k \eta_{a_i})$ are
defined by the recursion relation
\begin{equation}
(\sum_{i=1}^k a_i)(k-1) W_0(\prod_{i=1}^k \eta_{a_i})=
\sum_{k<l} (a_k+a_l) W_0(\eta_{a_k+a_l}\prod_{i\not= k,l}
 \eta_{a_i}),
\end{equation}
and initial data
\begin{equation}
W_0(\eta_{a_1}\eta_{a_2})=1
\end{equation}
(we give here the special case of this relation for genus $0$).

Note that the numbers $W_0(\prod_{i=1}^k \eta_{a_i})$ are
uniquely determined by this relation. Note also that
$W_0(\prod_{i=1}^k \eta_{a_i})=1$ obviously satisfy
this relation and the initial data. So, for any $a_1,\dots,a_k$
$W_0(\prod_{i=1}^k \eta_{a_i})=1$. In
particular, $W_0(\eta_1^n)=1$.

\section{The other Hodge integrals}

In this section we give an algorithm to calculate any Hodge
integral
\begin{equation}
\int_{\oM_{g,1}}\psi_1^{3g-2-i}\lambda_i.
\end{equation}

\subsection{Algorithm}
By $W^i_g(\eta_{a_1}\dots\eta_{a_n})$ denote the intersection
number
\begin{equation}
W^i_g(\eta_{a_1}\dots\eta_{a_n}):=\int_{V_g(a_1,\dots,a_n)}
\psi_1^{2g+n-2-i}\lambda_i
\end{equation}

The first step of the algorithm is just the same as in the
$\lambda_g$-case.

\begin{theorem}\label{binomialG}
For arbitrary positive integers
$a_1,\dots,a_n$ we have
\begin{multline}
(-1)^g g! \int_{\oM_{g,1}}\psi_1^{3g-2-i}\lambda_i =
\binom{g}{0}W_g^i(\prod_{i=1}^n \eta_{a_i})
- \binom{g}{1}
W_g^i(\eta_1\prod_{i=1}^n \eta_{a_i})
+ \\
\binom{g}{2}
W_g^i(\eta_1^2\prod_{i=1}^n \eta_{a_i})
-
\dots
+ (-1)^g \binom{g}{g}
W_g^i(\eta_1^g\prod_{i=1}^n \eta_{a_i})
\end{multline}
\end{theorem}

The recursion relation for the numbers
$W_g^i(\prod_{i=1}^n \eta_{a_i})$ looks like follows:

\begin{theorem}\label{recursionG}
If $2g+n-2-i>0$, then
\begin{multline}\label{recG}
(\sum_{i=1}^n a_i)(2g+n-1) W_g^i(\prod_{i=1}^n \eta_{a_i})=
\sum_{k<l} (a_k+a_l) W_g^i(\eta_{a_k+a_l}\prod_{i\not= k,l}
 \eta_{a_i})\\
+
\sum_{k=1}^n \frac{a_k^3-a_k}{12}
W_{g-1}^{i-1}(\prod_{i=1}^n \eta_{a_i})
+
\frac{1}{2}\sum_{k=1}^n
\sum_{a_k'+a_k''=a_k} a_k'a_k''
W_{g-1}^{i}(\eta_{a_k'}\eta_{a_k''}\prod_{i\not=k} \eta_{a_i}).
\end{multline}
\end{theorem}

The initial values are just the same as in the $\lambda_g$-case.

\begin{theorem}\label{initialG}
If $i>g$ then $W^i_g(\eta_{a_1}\dots\eta_{a_n})=0$. Besides,
\begin{equation}
W_1^1(\eta_{a_1})=\frac{a_1^2-1}{24}; \qquad
W_0^0(\eta_{a_1}\eta_{a_2})=1.
\end{equation}
\end{theorem}

Proofs of these Theorems are very similiar to the proofs of
Theorems~\ref{binomial}-\ref{initial}. So we skip the proofs.

\subsection{Some calculations}
Let us apply our algorithm to compute the integral
$\int_{\oM_{2,1}}\psi_1^3\lambda_1$. Using
Theorem~\ref{binomialG}, we get
\begin{equation}
\int_{\oM_{2,1}}\psi_1^3\lambda_1=\frac{1}{2}W^1_2(\eta_1^3)-
W^1_2(\eta_1^2).
\end{equation}

Applying many times Theorems~\ref{recursionG}
and~\ref{initialG}, we obtain
\begin{equation}
W^1_2(\eta_1^3) = \frac{1}{120};\qquad
W^1_2(\eta_1^2) = \frac{1}{480}.
\end{equation}

Therefore,
\begin{equation}
\int_{\oM_{2,1}}\psi_1^3\lambda_1=\frac{1}{240}-
\frac{1}{480}=\frac{1}{480}.\\
\end{equation}

The same number is given by Equation~\ref{FPformula}.

\end{document}